\documentclass[12pt]{amsart}

\usepackage{srcltx}  

\usepackage{geometry} 
\usepackage{amsmath}
\usepackage{amsfonts,amssymb,amsthm,amscd,latexsym,euscript}

\usepackage{xypic}
\xyoption{all}
\def\umono{\ar@{_{(}->}[u]}
\def\uumono{\ar@{_{(}->}[uu]}

\def\lmono{\ar@{_{(}->}[l]}
\def\llmono{\ar@{_{(}->}[ll]}

\newcommand{\map}{\operatorname{map}\nolimits}
\newcommand{\op}{{\mathrm{op}}}

\newcommand{\A}{\ifmmode{\mathcal{A}}\else${\mathcal{A}}$\fi}
\newcommand{\K}{\ifmmode{\mathcal{K}}\else${\mathcal{K}}$\fi}
\newcommand{\U}{\ifmmode{\mathcal{U}}\else${\mathcal{U}}$\fi}
\newcommand{\T}{\ifmmode{\mathcal{T}}\else${\mathcal{T}}$\fi}
\newcommand{\FF}{\ifmmode{\mathcal{F}}\else${\mathcal{F}}$\fi}
\newcommand{\LL}{\ifmmode{\mathcal{L}}\else${\mathcal{L}}$\fi}

\newtheorem{theorem}{Theorem}[section]
\newtheorem{proposition}[theorem]{Proposition}
\newtheorem{corollary}[theorem]{Corollary}
\newtheorem{lemma}[theorem]{Lemma}

\theoremstyle{definition}
\newtheorem{definition}[theorem]{Definition}
\newtheorem{example}[theorem]{Example}
\theoremstyle{remark}
\newtheorem{remark}[theorem]{Remark}

\def\N{{\mathbb{N}}}

\DeclareMathOperator{\Nil}{{\textsl{Nil}}}
\DeclareMathOperator{\Ho}{{\textup{Ho}}}

\newcommand{\sS}{{sSets}}

\newcommand{\cat}[1]{\ensuremath{\EuScript #1}}

\title[Goodwillie calculus and Whitehead products]{Goodwillie calculus and Whitehead products}

\author{Boris Chorny and J\'er\^{o}me Scherer}

\thanks{The second author is supported by FEDER/MEC grant MTM2010-20692.}

\thanks{}


\begin{document}


\begin{abstract}
We prove that iterated Whitehead products of
length $(n+1)$ vanish in any value of an $n$-excisive functor
in the sense of Goodwillie. We compare then different notions
of homotopy nilpotency, from the Berstein-Ganea definition to
the Biedermann-Dwyer one. The latter is strongly related to
Goodwillie calculus and we analyze the vanishing of iterated 
Whitehead products in such objects.
\end{abstract}


\maketitle


\section*{Introduction}
\label{sec:intro}
Goodwillie calculus, \cite{MR1162445}, \cite{MR2026544}, gives a systematic
way to approximate a functor (say from spaces to spaces) by a tower of
functors satisfying higher excision properties. When applied to the identity functor
this tower reflects remarkable periodicity properties, as investigated by
Arone and Mahowald, \cite{MR1669268}. More recently Biedermann
and Dwyer, \cite{MR2580428}, used the stages of the very same tower to construct 
(simplicial) algebraic theories in the sense of Lawvere, \cite{MR0158921}.
The homotopy algebras over these theories are called homotopy
nilpotent groups, and the class of nilpotency corresponds exactly to the
chosen stage of the Goodwillie tower.

Our objectives in this article are twofold. First we investigate why $n$-excisive
functors should be related to homotopy nilpotency in the classical sense.
In the early sixties, Berstein and Ganea introduced a concept of nilpotent
loop spaces, \cite{MR0126277}. They require that an iterated commutator
map be trivial up to homotopy, which implies in particular that iterated Samelson
products vanish in the loop space $\Omega X$, or equivalently, that iterated 
Whitehead products vanish in~$X$. Already G.~Whitehead, \cite{MR0016672}, 
had the insight that the (J.H.C.) Whitehead products satisfy identities which reflect 
commutator identities for groups.  Work of Hopkins, \cite{MR1017164},
drew renewed attention to such questions by relating this classical nilpotency notion with
the nilpotence theorem of Devinatz, Hopkins, and Smith, \cite{MR960945}.
We prove the following.

\medskip

\noindent {\bf Theorem~\ref{prop Whitehead products and excision}}
{\it Let $F$ be any $n$-excisive functor from the category of 
pointed spaces to pointed spaces. Then all $(n+1)$-fold iterated
Whitehead products vanish in~$F(X)$ for every finite space~$X$.
}

\medskip

Our result shows in fact that $\Omega F(X)$ is a homotopy nilpotent loop space in the sense of Ganea and 
Bernstein for every $n$-excisive functor $F$ and every finite space~$X$. 

The difficulty of the proof resides in finding a way to take into account the global
property of the functor (to be $n$-excisive) and not to focus on a particular
value $F(X)$. Except for this, the proof uses the general theory of
Goodwillie calculus.

In the second part of the article we look more closely at the relationship
between the different types of homotopy nilpotency available on the market.
We start with the classical algebraic theory $\Nil_n$ describing nilpotent groups
of class $\leq n$, and observe that Berstein-Ganea nilpotent loop spaces are 
$\Nil_n$-algebras in the homotopy category of spaces. We show that homotopy
$\Nil_n$-algebras in the sense of Badzioch, \cite{MR1923968}, are always
homotopy nilpotent in the sense of Biedermann and Dwyer. Finally, both are
$\Nil_n$-algebras in the homotopy category of spaces, so that the following theorem 
applies to all kinds of homotopy nilpotent groups that appeared so far in the literature,
and in particular to the Biedermann-Dwyer ones.

\medskip

\noindent {\bf Theorem~\ref{theorem Whitehead products}}
{\it Let $\Omega X$ be a homotopy nilpotent group of class $\leq n$.
Then all $(n+1)$-fold iterated Whitehead products vanish in~$X$.}

\medskip

The proof depends on a non-trivial computation of sets of components 
in \cite{MR2580428}. As a corollary, since values of $n$-excisive functors
yield examples of homotopy nilpotent groups of class $\leq n$, we obtain
another proof of Theorem~\ref{prop Whitehead products and excision}.

\medskip

\noindent {\bf Acknowledgments.} We would like to thank Amnon Neeman and the Mathematical Sciences Institute
at ANU where this project started. The hospitality of EPFL is greatly acknowledged by the first author, who visited 
the second author in order to complete this work. We learned about homotopy nilpotent groups from Georg Biedermann and Bill
Dwyer and thank them for their interest in this project.

\section{Samelson and Whitehead products}
\label{sec:products}
We recall briefly the definition of Samelson and Whitehead products and construct
a ``universal space" built from wedges of spheres in which higher Whitehead products
vanish.

Let $X$ be a pointed space. Given $\alpha \in \pi_{a+1} X$ and
$\beta \in \pi_{b+1} X$, take the adjoint classes $\alpha' \in \pi_a
(\Omega X)$ and $\beta' \in \pi_b (\Omega X)$. The composite of the
product map $\alpha' \times \beta': S^a \times S^b \rightarrow
\Omega X \times \Omega X$ with the commutator map $\Omega X \times
\Omega X \rightarrow \Omega X$ is null-homotopic when restricted to
the wedge $S^a \vee S^b$ and thus factors through $S^{a+b}$, uniquely
up to homotopy. This
factorization represents the Samelson product $\langle \alpha', \beta'
\rangle \in \pi_{a+b} \Omega X$ and the adjoint class is the Whitehead 
product $[\alpha, \beta] \in \pi_{a+b+1} X$. 

\begin{remark}
\label{remark order}
{\rm Iterated Whitehead  products can be computed as adjoint to iterated Samelson products. For example
a triple Whitehead product of the form $[[\alpha, \beta], \gamma]$ coincides with
the adjoint of the Samelson product $\langle \, \langle \alpha', \beta' \rangle, \gamma' \rangle$.
Let us also mention that the order of the classes in a Whithead product does not matter (up to a sign). 
We will therefore concentrate on one standard choice of bracketing.}
\end{remark}

By definition, the Whitehead product $[\iota_1, \iota_2]$ of the two canonical inclusions
$\iota_1: S^a \hookrightarrow S^a \vee S^b$ and $\iota_2: S^b \hookrightarrow S^a \vee S^b$
is the attaching map of the top cell in $S^a \times S^b$. Moreover, any Whitehead product
$[\alpha, \beta]: S^{a+b+1} \rightarrow X$ factors through $[\iota_1, \iota_2]$. This motivates
the construction of a space built from wedges of spheres which will be crucial for
understanding when certain iterated Whitehead products vanish. 
We consider  $n+1$ positive integers $k_1, \dots , k_{n+1}$ and the
wedge of $n+1$ spheres $W= \vee S^{k_i}$. Denote by $\iota_i:
S^{k_i} \rightarrow W$ the wedge summand inclusion. Define the
$(n+1)$-cube of pointed spaces $V: \mathcal P(\underline{n+1})\setminus \{\emptyset\}
\rightarrow Spaces_*$ by sending a subset $S \subset
\underline{n+1}$ to $\bigvee_{i\not\in S} S^{k_i}$. Adding $W$ as initial value 
$V(\emptyset)$ makes this diagram a strongly homotopy co-Cartesian cube
(we will also write abusively $V(i)$ instead of $V(\{i\})$ to ease the notation). 
We let $Q$ be the homotopy inverse limit of $V$, and to fix a representative
we take $Q$ to be the inverse limit of the fibrant replacement 
$V \tilde\hookrightarrow \hat V$ of this diagram in the injective model structure, \cite{Heller, Jardine}.

\begin{example}
\label{example inverse}
{\rm When $n=1$, we have two spheres $S^{k_1}$ and $S^{k_2}$. The diagram $V$ is the pull-back diagram
$S^{k_1} \rightarrow * \leftarrow S^{k_2}$ and $Q=S^{k_1} \times S^{k_2}$. The Whitehead product of the
summand inclusions is trivial in $Q$.}
\end{example}

The looped diagram $\Omega V$ is easier to analyze since the loop space on a 
wedge of spheres splits by the Hilton-Milnor theorem, see the original article
\cite{MR0068218} or Milnor's unpublished article in
\cite{MR0445484}: Each ``basic word" $w$ in $x_1, \dots, x_{n+1}$ determines
a Whitehead product in $\pi_{N(w)} (S^{k_1} \vee \dots \vee S^{k_{n+1}})$
and $\Omega (S^{k_1} \vee \dots \vee S^{k_{n+1}}) \simeq \prod_w \Omega S^{N(w)}$.
Thus, when $n=2$, the basic word $x_1x_2x_3$ corresponds to the Whitehead product
$[[\iota_1, \iota_2], \iota_3]$ represented by a map $S^{k_1 + k_2 + k_3 -2} \rightarrow S^{k_1} \vee S^{k_2}  \vee S^{k_{3}}$.

\begin{lemma}
\label{lemma wedges of spheres}
The loop space $\Omega Q$ is homotopy equivalent to a product of loop spaces on
spheres, namely $\displaystyle \prod \Omega S^{N(w)}$ where the product is taken
over all basic words in at most $n$ of the letters $x_1, \dots, x_{n+1}$.
\end{lemma}

\begin{proof}
We identify $\Omega Q$ with the homotopy inverse limit of the diagram $\Omega V$,
each value of which splits as a product of loop spaces on spheres:
\[
\Omega V(S) \simeq \prod_{i \not\in S} \Omega S^{k_i} \times \dots \times \prod_{w \in W_S} \Omega S^{N(w)}
\]
where $W_S$ is the subset of those basic words written with all $x_i$'s with $i\not\in S$. 
We observe that each map $\Omega V(S) \rightarrow \Omega V(T)$, with $S \subset T$, is the projection
on the summands $\Omega S^{N(w)}$ corresponding to the basic words not written
with the letters in $T$. Therefore the diagram $\Omega V$ is a hypercube of which the
homotopy inverse limit is the product of all $\prod_{w \in W_S} \Omega S^{N(w)}$ with 
$S \neq \emptyset$.
\end{proof}

For any choice of bracketing $n+1$ elements there is an
$(n+1)$-fold Whitehead product $w: S^{k_1 + \dots + k_{n+1} -n}
\rightarrow W$. We denote by $C_w$ the homotopy cofiber 
of $w$.

\begin{lemma}
\label{lemma vanishing}
The $(n+1)$-fold Whitehead product $[ [ \dots [ [\iota_1, \iota_2], \iota_3], \dots], \iota_{n+1}]$ vanish in $Q$. 
\end{lemma}

\begin{proof}
This Whitehead products vanish in $Q$ if and only if
the adjoint Samelson product  vanish in $\Omega Q$.
Since $\Omega Q$ splits as a product of loop spaces on spheres, it is
sufficient to prove that the projection on each factor is
null-homotopic. By Lemma~\ref{lemma wedges of spheres} each 
factor already appears in $\Omega V(S)$ for some non-empty subset $S$,
so that, by adjunction again, it is enough to show that the image in $V(i)$ 
of our $(n+1)$-fold Whitehead product  vanishes for any $1\leq i \leq n+1$. This is
so because the image of $\iota_i$ in $V(i)$ is the trivial map and any
Whitehead product involving the trivial map is null-homotopic.
\end{proof}




\section{The values of $n$-excisive functors}
\label{sec:excision}
We perform our main computation in this section. Let $F$ be an $n$-excisive
functor from pointed spaces to pointed spaces (so $F$ sends strongly homotopy
co-Cartesian $(n+1)$-cubes to homotopy Cartesian ones. We prove that
all $(n+1)$-fold Whitehead products vanish in $F(X)$ for any space $X$.
Because it is very difficult to use the global property of excision by focusing
on one single value of the functor $F$, 
we will use pointed representable functors $R^X$, defined by $R^X(Y) = \map_*(X, Y)$. 
For any pointed space $A$, a natural transformation $R^X \wedge A \rightarrow F$ 
corresponds by adjunction to a map $A \rightarrow \hbox{\rm hom}(R^X, F)$, i.e. to a
map $A \rightarrow F(X)$ by the enriched Yoneda Lemma \cite[2.31]{Kelly}.

\begin{theorem}
\label{prop Whitehead products and excision}
Let $F$ be any $n$-excisive functor from the category of
pointed spaces to pointed spaces. Then all $(n+1)$ fold iterated
Whitehead products vanish in~$F(X)$ for every finite space $X$.
\end{theorem}

\begin{proof}
Let us fix homotopy classes of maps $\alpha_i: S^{k_i} \rightarrow F(X)$ for
$1 \leq i \leq n+1$. We need to prove that the iterated Whitehead product
$[ [ \dots [[\alpha_1, \alpha_2], \alpha_3], \dots ], \alpha_{n+1}] $ is zero. This
product is represented by a map
\[
S^{k_1 + \dots k_{n+1} + 1} \xrightarrow{w} \vee_{i=1}^{n+1} S^{k_i} = W \rightarrow F(X)
\] 
which is null-homotopic if it factors through the homotopy cofiber $C_w$ of the
``universal" $(n+1)$-fold Whitehead product $w$. The use of representable functors
translates then as follows: We need to show that
any natural transformation $\eta: R^X \wedge W \rightarrow F$ factors
through $R^X \wedge C_w$. As $F$ is $n$-excisive,
there exists a natural transfortmation $P_n(R^X \wedge W) \rightarrow F$ 
such that the composite $R^X \wedge W \rightarrow P_n(R^X \wedge W) \rightarrow F$
coincides with $\eta$ up to homotopy.
It is thus enough to construct a natural transformation $R^X
\wedge C_w \rightarrow P_n(R^X \wedge W)$.

Smashing the diagram $V$ with a representable functor we obtain a hypercube
$R^X \wedge V$ of functors, which is strongly homotopy cocartesian since $V$
is so. We focus on the natural transformations $R^X \wedge W \rightarrow R^X \wedge V(i)$.
If $c = \hbox{\rm dim} X$, and $Y$ is a $k$-connected space with $k \geq
c$, then $R^X(Y)$ is $(k-c)$-connected and $(R^X \wedge W)(Y) \rightarrow (R^X \wedge V(i))(Y)$
is $(k - c + k_i)$-connected.
Let $G$ denote the homotopy inverse limit of the diagram of functors $R^X \wedge V$.

The generalized Blackers-Massey theorem
\cite[Theorem~2.3]{MR1162445} implies that the natural
transformation $\theta: R^X \wedge W \rightarrow G$ is $[(n+1)k -(n+1)c + \sum k_i -n]$-connected
when evaluated at a $k$-connected space with $k \geq c$. This implies that $R^X \wedge
W$ and $G$ agree to order $n$ in the terminology of
\cite[Definition~1.2, Proposition~1.6]{MR2026544},
 so that $P_n(R^X \wedge W) \simeq
P_n(G)$.

Lemma~\ref{lemma vanishing} yields a map $C_w \rightarrow Q$
such that $W \rightarrow C_w \rightarrow Q$ is the natural map from
$W$ to the homotopy inverse limit of the diagram $V$ (we fix the model
$C_w = W \cup_w D^{k_1 + \dots + k_{n+1} +2}$ for the homotopy
cofiber so that the factorization is strict). We interpret this map as
 a map from the constant diagram $C_w$ to a fibrant replacement $\hat V$ of~$V$
 in the injective model category of hypercubical diagrams.
Smashing with a representable functor we get a natural
transformation $R^X \wedge C_w \rightarrow R^X \wedge V$.
Taking homotopy inverse limits we obtain finally a natural transformation
$R^X \wedge C_w \rightarrow G$ such that the the composite
$R^X \wedge W \rightarrow R^X \wedge C_w \rightarrow G$
coincides with $\theta$. The natural transformation
\[
R^X \wedge C_w \rightarrow G \rightarrow P_n G \simeq P_n(R^X \wedge W)
\]
is the one we needed to conclude.
\end{proof}

\begin{remark}
\label{rem generalize}
{\rm This proof easily generalizes to show that iterated
\emph{generalized} Whitehead products vanish. It suffices to replace the
Hilton splitting theorem for loop spaces on a wedge of spheres by
Milnor's generalized version for wedges of suspensions.}
\end{remark}

\section{Nilpotent groups and algebraic theories}
\label{sec:groups}
Let us first recall the classical concept of algebraic theory due to Lawvere \cite{MR0158921} and some of its modern variations.

\begin{definition}
\label{def theory}
A small category $T$ is an algebraic theory if the objects of $T$ are indexed by natural numbers 
$\{T_{0}, T_{1},\ldots, T_{n},\ldots\}$ and for all $n\in \N$ the $n$-fold categorical coproduct of $T_{1}$ is 
naturally isomorphic to $T_{n}$. The algebraic theory $T$ is \emph{simplicial} if it is a (pointed) simplicial 
category, i.e., $T$ is enriched over $\sS_{\ast}$

Let $\cat C$ be  a category. A $\cat C$-\emph{algebra} over a theory $T$ is a functor $A\colon T^{\op}\to \cat C$ 
taking coproducts in $T$ into products in $\cat C$. 

If $T$ is a simplicial algebraic theory and $\cat C=\sS_{\ast}$, then we distinguish between \emph{strict} and 
\emph{homotopy} simplicial algebras, which are simplicial functors $A\colon T^{\op}\to \cat C$ taking 
coproducts in $T$ to products in \cat C strictly or up to homotopy, respectively.
\end{definition}

The categories of simplicial algebras and homotopy simplicial algebras were 
compared by Badzioch in \cite{MR1923968}.
He proved that any homotopy algebra can be rigidified to a strict algebra.

Of central interest for us will be algebras over algebraic theories defined in the 
homotopy category of simplicial  sets $\cat C=\Ho(\sS_{\ast})$. We call them \emph{algebras up to homotopy}, 
in order to distinguish them from the homotopy algebras defined above. 
There is a natural way to associate to every homotopy algebra $A$, an algebra up to homotopy: just compose the 
functor $A$ with the product preserving functor $\Gamma\colon \sS_{\ast}\to \Ho(\sS_{\ast})$. Formally,
we need to choose homotopy inverse maps $f_k: A(k) \rightarrow A(1)^k$ and $g_k: A(1)^k \rightarrow A(k)$
and replace each morphism $A(h): A(m) \rightarrow A(n)$ by the composite $f_n \circ A(h) \circ g_m$.
The converse is not true of course, and we will encounter examples of algebras up to homotopy
which cannot be upgraded to homotopy algebras. 

Lawvere in his seminal article \cite{MR0158921} has discovered the fundamental fact that an algebraic theory defining
a variety as the category of algebras, is the dual of the subcategory of finitely generated free algebras in this variety.
In this work we will look closer into the algebraic theories defining the concepts of groups and nilpotent groups of 
class $\leq n$ in various settings. 

Consider thus the full subcategory $\Nil_n$ of the category of groups: the objects
are the free nilpotent groups $F_k/ \Gamma_{n+1} F_k$ of class $n$,
for all $k\geq 1$. In the category of nilpotent groups of class
$\leq n$, these  groups are free in the sense that they can be identified
with the coproducts of $k$ copies of $\mathbf Z = F_1/\Gamma_{n+1} F_1$.  
The set of morphisms from $1$ to $k$ is precisely the group $F_k/
\Gamma_{n+1} F_k$. When $n=\infty$, we define the objects of $\Nil_\infty$
to be the free groups $F_k$. A $\Nil_n$-\emph{algebra} in $Sets$ is thus
a product preserving contravariant functor $N\colon \Nil_n^{\op} \rightarrow Sets$.

\begin{proposition}
\label{prop nilpotent group}
A $\Nil_n$-algebra is a nilpotent group of class $\leq n$.
\end{proposition}

Because it will play an important role in the sequel, let us be
precise and say explicitly how the group structure arises and why it is nilpotent.
By abuse of notation we write also $N$ for the value $N(1)$.
The multiplication $m: N \times N \rightarrow N$ is the morphism
corresponding to the product of the two generators of $F_2$ in the
quotient $F_2/ \Gamma_{n+1} F_2$ and the inverse is the morphism $N
\rightarrow N$ corresponding to the inverse of the generator of
$F_1$. It is easy to check that this equips $N$ with a group
structure. This is in fact equivalent to the structure of a $\Nil_\infty$-algebra:
Given $k$ elements $n_1, \dots, n_k \in N$ and a word $w$ in $k$ letters,
the product $w(n_1, \dots, n_k)$ can be read of from the
morphism $N^k \rightarrow N$ corresponding to $w$.
The claim about the nilpotency class follows then from
the fact that all words of the form $[ [\dots [[x_1, x_2], x_3], \dots ], x_{n+1}]$
are identified to $1$ in $F_{n+1}/ \Gamma_{n+1} F_{n+1}$. Hence any iterated commutator of
length $\geq n+1$ must be trivial in a $Nil_n$-algebra.

\begin{remark}
\label{rem simplicial nilpotent group}
{\rm A $\Nil_n$-algebra in the category of simplicial sets, i.e. a
product preserving contravariant functor $N\colon \Nil_n^{op} \rightarrow
sSets$, is a simplicial nilpotent group of class $\leq n$. In
particular when $n=1$ we are considering simplicial abelian groups,
i.e. generalized Eilenberg-Mac Lanes spaces, so called ``GEMs", see
for example \cite{MR0279808}. Schwede also considers such
objects and compares them stably, \cite[Example~7.4]{MR1791267}
with a category of modules over a Gamma-ring. 

Badzioch's rigidification result states in this context that any homotopy $\Nil_n$-algebra is
homotopy equivalent to a strict $\Nil_n$-algebra. Again for $n=1$, this means 
that all homotopy $\Nil_1$-algebras are homotopy equivalent to GEMs.
This is not quite what we would like
to study when we are speaking about a homotopy version of abelian
topological groups (what we understand under this name is rather an
infinite loop space). The notion of $\Nil_n$-algebras in simplicial sets
is thus too rigid and we will need to relax it a little.
}
\end{remark}

\section{Nilpotent groups in the homotopy category}
\label{sec:homotopy}
In the next section we will turn to the solution Biedermann and Dwyer
found to describe homotopy nilpotency. But before we do that, we first
describe the most naive way to define nilpotency in homotopy theory.

\begin{definition}
\label{def nilpotent space}
{\rm A \emph{nilpotent group up to homotopy of class $\leq n$} is a
product preserving contravariant functor $N: \Nil_n^{op} \rightarrow
\hbox{\rm Ho}(Spaces_*)$.}
\end{definition}

How do these nilpotent groups up to homotopy look like? They are
pointed spaces $G$ together with a homotopy associative
multiplication and a homotopy inverse (i.e. group-like $H$-spaces)
coming from the morphisms in  $\Nil_n^{op}(2, 1)$ and $\Nil_n^{op}(1,
1)$ described in the previous section, such that all higher
commutator maps of order $n+1$ are null-homotopic. Berstein and Ganea, 
\cite[Definition~1.7]{MR0126277} give a definition of nilpotency for group like 
spaces by requiring that the
$(n+1)$-st commutator map be null-homotopic. Their work predates by
two years the introduction by Lawvere of algebraic theories, and is
therefore not stated in the language we have used, but it is equivalent.

\begin{proposition}
\label{def nilpotent group}
A nilpotent group up to homotopy is a homotopy nilpotent
group in the sense of Berstein and Ganea. \hfill{\qed}
\end{proposition}

\begin{example}
\label{exm abelian up to homotopy}
{\rm When $n=1$, a loop space is abelian (nilpotent of class $\leq 1$) up to homotopy
if the commutator map $\Omega X \times \Omega X \rightarrow \Omega X$ is null-homotopic,
i.e. if the product is homotopy commutative. Thus any double loop space is abelian up
homotopy. When $n=\infty$, groups up to homotopy are simply group objects in the
homotopy category, i.e. homotopy associative $H$-spaces with inverse.}
\end{example}

These examples show that the Berstein-Ganea definition is too flexible. When looking
at loop spaces, the filtration given by nilpotency up to homotopy interpolates roughly between
loop spaces and double loop spaces. However it allows us to read off the vanishing of
iterated Samelson products. This is basically \cite[Theorem~4.6]{MR0126277}.

\begin{proposition}
\label{prop vanishing of Whitehead products}
Let $X$ be a pointed space and assume that the loop space $\Omega X$
is nilpotent up to homotopy of class $\leq n$. Then all $(n+1)$ fold
iterated Whitehead products vanish in~$X$.
\end{proposition}

\begin{proof}
The vanishing of iterated Whitehead products is equivalent to the
vanishing of iterated Samelson products in the loop space. This
follows now directly from the fact that in a $\Nil_n$-algebra in the
homotopy category the $(n+1)$-fold commutator map $(\Omega X)^{n+1}
\rightarrow \Omega X$ is null-homotopic by definition.
\end{proof}

\begin{example}
\label{exm Porter}
{\rm Porter proved that $S^3$ is nilpotent up to homotopy of class $3$, \cite{MR0169244}. 
There is a non-trivial three fold Whitehead product in $BS^3$, but all $4$-fold products vanish.
However, the compact Lie group $S^3$ is not nilpotent as a group. More generally, Rao
proved that compact Lie groups are nilpotent up homotopy if and only if they are torsion free, \cite{MR1441487}.
The if part is due to Hopkins, \cite{MR1017164}.
}
\end{example}

\section{Enriched homotopy nilpotent groups}
\label{sec:biedermann and dwyer}
This section finally introduces the ``correct" homotopy nilpotent groups.
We recall their definition, show that iterated Samelson products vanish
in such spaces, and compare them to homotopy $\Nil_n$-algebras and
spaces which are nilpotent up to homotopy in the sense of Berstein and Ganea.

In their recent work \cite{MR2580428} Biedermann and Dwyer define homotopy
nilpotent groups as homotopy $\mathcal G_n$-algebras in the category
of pointed spaces, where $\mathcal G_n$ is a \emph{simplicial} algebraic
theory constructed from the Goodwillie tower of the identity.
Concretely, the object $k$ corresponds to the $k$ fold product of
the functor $\Omega P_n(id)$ in the category of functors from finite
pointed spaces to pointed spaces. Hence a homotopy nilpotent group
of class $\leq n$ is the value at $1$ of a simplicial functor
$\tilde X$ from $\mathcal G_n$ to pointed spaces which commutes up
to homotopy with products. Homotopy algebras in an enriched
context have been studied by Rosick{\'y} in~\cite{MR2349711}.

\begin{proposition}
\label{prop comparison}
A homotopy $\Nil_n$-algebra is always a homotopy $\mathcal G_n$-algebra.
Both of them are $\Nil_n$-algebras up to homotopy.
\end{proposition}

\begin{proof}
The space of maps from $k$ to $1$ in the algebraic theory $\mathcal G_n$, 
which is by definition the space of all natural transformations from 
$(\Omega P_n(id))^k$ to $\Omega P_n(id)$, is identified as
the space $\Omega P_n(id)(\vee_k S^1)$, \cite[Corollary~4.7]{MR2580428}.
Biedermann and Dwyer's main computation shows that
the group of connected components coincides
with the free nilpotent group of class $n$ on $k$ generators:
$$
\pi_0 \mathcal P_n(k, 1) \cong \pi_0 \bigl(\Omega P_n(id)(\vee_k
S^1)\bigr) \cong F_k/ \Gamma_{n+1} F_k.
$$
There is hence  a functor of simplicial algebraic theories $\pi_0: \mathcal G_n \rightarrow \Nil_n$.
Thus any homotopy $\Nil_n$-algebra can be seen as a homotopy $\mathcal G_n$-algebra
by pulling back along $\pi_0$.

Consider now a homotopy nilpotent group $\Omega X$ of class $\leq n$ given
as the value at $1$ of a homotopy $\mathcal G_n$-algebra 
$\tilde X: \mathcal G_n \rightarrow Spaces_*$. The composite
diagram $F: \mathcal G_n \rightarrow Spaces_* \rightarrow \hbox{\rm
Ho}(Spaces_*)$ is now simply a diagram $F: \Nil_n \rightarrow \hbox{\rm
Ho}(Spaces_*)$ as we keep from the simplicial data only
one homotopy class of maps $\tilde X(k) \rightarrow \tilde X(l)$ for
each connected component of ${\mathcal G_n}(k, l) \simeq \Omega
P_n(id)(\vee_k S^1)^l$. The second claim then follows from the general
procedure we described in Section~\ref{sec:groups} to get an algebra
up to homotopy from a homotopy algebra.
\end{proof}

\begin{theorem}
\label{theorem Whitehead products}
Let $\Omega X$ be a homotopy nilpotent group of class $\leq n$. Then
all $(n+1)$ fold iterated Whitehead products vanish in~$X$.
\end{theorem}

\begin{proof}
The Berstein-Ganea Proposition~\ref{prop vanishing
of Whitehead products} implies the vanishing of all iterated
Whitehead products of length $n+1$ in $X$.
\end{proof}

\begin{remark}
\label{remark loop structure}
{\rm Observe here that a homotopy nilpotent group of class $n$ is
also a homotopy nilpotent group of class $\infty$ since we have a
map of algebraic theories $\mathcal P_\infty \rightarrow \mathcal
P_n$. This means that a homotopy nilpotent group of class $n$ has
the homotopy type of a loop space and the multiplication derived
from the algebraic theory is this precise loop multiplication. This
is what allows us to use the Berstein-Ganea result in the last line
of the previous proof.}
\end{remark}

\begin{example}
\label{example BD}
{\rm Homotopy abelian groups, that is homotopy $\mathcal G_1$-algebras, are infinite
loop spaces and homotopy groups, i.e. homotopy $\mathcal G_\infty$-algebras, are
loop spaces. This is why the notion of homotopy nilpotency of Biedermann and Dwyer
is better than the others. It interpolates between the ``right" notions of groups and abelian 
groups in homotopy theory. In particular, $BU$ is
homotopy abelian, but not a homotopy $\Nil_1$-algebra, and $\Omega^2 S^4$
is abelian up to homotopy, but not a homotopy abelian group. This
illustrates how the different notions of nilpotency differ.
}
\end{example}

\begin{corollary}
\label{corollary Whitehead products and excision}
Let $F$ be any $n$-excisive functor from the category of
pointed spaces to pointed spaces. Then all $(n+1)$ fold iterated
Whitehead products vanish in~$F(X)$ for any finite space~$X$.
\end{corollary}

\begin{proof}
Biedermann and Dwyer prove that $n$-excisive functors produce
examples of homotopy nilpotent groups: $\Omega F(X)$ is homotopy nilpotent
of class $\leq n$, \cite[Corollary~9.3]{MR2580428}.
\end{proof}

\begin{remark}
\label{rem iff}
{\rm Biedermann and Dwyer claim after \cite[Corollary~9.3]{MR2580428}
that all homotopy nilpotent groups are given as values of loops on $n$-excisive 
functors. Combined with Theorem~\ref{prop Whitehead products and excision}
this gives another proof of the fact that homotopy
nilpotent group of class $n$ have vanishing $(n+1)$-fold iterated
Whitehead products.}
\end{remark}


\bibliographystyle{amsplain}
\providecommand{\bysame}{\leavevmode\hbox to3em{\hrulefill}\thinspace}
\providecommand{\MR}{\relax\ifhmode\unskip\space\fi MR }
\providecommand{\MRhref}[2]{%
  \href{http://www.ams.org/mathscinet-getitem?mr=#1}{#2}
}
\providecommand{\href}[2]{#2}



\bigskip
{\small
\begin{minipage}[t]{8 cm}
Boris Chorny\\ 
Department of Mathematics\\
University of Haifa at Oranim\\
IL - 36006 Qiryat Tivon, Israel\\
\textit{E-mail:}\texttt{\,chorny@math.haifa.ac.il}\\
\end{minipage}
\begin{minipage}[t]{8 cm}
J\'er\^ome Scherer\\
EPFL SB MATHGEOM\\
Station 8, MA B3 455\\
CH - 1015 Lausanne, Switzerland\\
\textit{E-mail:}\texttt{\,jerome.scherer@epfl.ch}\\
\end{minipage}

\end{document}